# ON HERMITE-HADAMARD TYPE INEQUALITIES VIA $h$-CONVEXITY

Muhammad IQBAL[1], Muhammad MUDDASSAR[2*] and Muhammad IQBAL BHATTI[3]

[1]Assistant Professor, Department of Mathematics, Govt. College of Technology, Sahiwal, Pakistan (miqbal.bki@gmail.com)

[2]Assistant Professor, Department of Mathematics, University of Engineering & Technology, Taxila, Pakistan (malik.muddassar@gmail.com)

[3]Professor, Department of Mathematics, University of Engineering & Technology, Lahore, Pakistan (uetzone@hotmail.com)

* (For Correspondence)

**ABSTRACT:** *Some new Hermite-Hadamard's inequalities for h-convex functions are proved, generalizing and unifying a number of known results. Some new applications for special Means of real numbers are also derived.*

## 1. Introduction

If a function $f : [a,b] \to \mathbb{R}$ is convex, induced

$$f\left(\frac{a+b}{2}\right) \leq \frac{1}{b-a} \int_a^b f(x)\,dx \leq \frac{f(a)+f(b)}{2} \qquad (1)$$

is known as the Hermite-Hadamard inequality. For concave function, $f$, the above order is reversed. Inequality (1) is refined, extended, generalized and new proofs are given in [1-11].

A midpoint type inequality in [10] is pointed out as follows:

**Theorem 1**. Let $f: I = [a,b] \neq \phi \subset \mathbb{R} \to \mathbb{R}$ twice differentiable function on $I^o$ with $f'' \in L^1(I)$. If $|f''|$ is convex function on $I$, induces

$$\left|\frac{f(a)+f(b)}{2} - \frac{1}{b-a}\int_a^b f(x)dx\right| \leq \frac{(b-a)^2}{24}\left[\frac{|f''(a)|+|f''(b)|}{2}\right] \quad (2)$$

In [1], A. Barani et al. proved a version for powers of a function whose absolute value of the second derivative is $P$-convex, for the right-hand side of (1) is incorporated in the following theorem as:

**Theorem 2.** Let $f: I = [a,b] \neq \phi \subset \mathbb{R} \to \mathbb{R}$ twice differentiable function on $I^o$ with $a < b$ and $f'' \in L^1(I)$. Assuming that $q \geq 1$ such that $|f''|^q$ is $P$-convex function on $I$, induces

$$\left|\frac{f(a)+f(b)}{2} - \frac{1}{b-a}\int_a^b f(x)dx\right| \leq$$

$$\frac{(b-a)^2}{24}\left[\left(|f''(a)|^q + \left|f''\left(\frac{a+b}{2}\right)\right|^q\right)^{\frac{1}{q}} + \left(|f''(b)|^q + \left|f''\left(\frac{a+b}{2}\right)\right|^q\right)^{\frac{1}{q}}\right] \quad (3)$$

A refinement of the left hand side of the Hermite-Hadamard inequality (1) is embodied in the following Theorem.

**Theorem 3**[2]. Let $f: I = [a,b] \neq \phi \subseteq \mathbb{R} \to \mathbb{R}$ be a convex function on $I$ with $a < b$, induces

$$\frac{1}{b-a}\int_a^b f(x)dx - f\left(\frac{a+b}{2}\right)$$
$$\geq \left|\frac{1}{b-a}\int_a^b \left|\frac{f(x)+f(a+b-x)}{2}\right|dx - \left|f\left(\frac{a+b}{2}\right)\right|\right|. \quad (4)$$

The following result gives us a refinement of second part of the Hermite-Hadamard inequality (1) for modulus.

**Theorem 4**[2]. Let $I$ be an interval of real numbers and $a, b \in I$ with $a < b$. If $f: I \subseteq \mathbb{R} \to \mathbb{R}$ *is* a convex function on $I$, then

$$\frac{f(a)+f(b)}{2} - \frac{1}{b-a}\int_a^b f(x)dx$$
$$\geq \begin{cases} |f(a)| - \frac{1}{b-a}\int_a^b f(x)dx, & \text{if } f(a) = f(b) \\ \frac{1}{f(b)-f(a)}\int_{f(a)}^{f(b)} |x|dx - \frac{1}{b-a}\int_a^b |f(x)|dx, & \text{if } f(a) \neq f(b) \end{cases}. \quad (5)$$

The following definitions are needed to explain the idea of $h$-convex function:

**Definition 1[2]:** Let $I$ be an interval of real numbers. A function $f: I \to \mathbb{R}$ is called convex if for all $x, y \in I$ and $t \in [0, 1]$, induces
$$f((1-t)x + ty) \leq (1-t)f(x) + tf(y)$$

$f$ is said to be concave, if the above inequality is reversed.

**Definition 2[5]:** A non $-$ negetive function $f: I \subseteq \mathbb{R} \to \mathbb{R}$ is called Godunova $-$ Levin function (or that $f \in Q(I)$) if for all $x, y \in I$ and $t \in (0, 1)$, implies
$$f(tx + (1-t)y) \leq \frac{1}{t}f(x) + \frac{1}{1-t}f(y)$$

It may be noted that non-negative monotone and non-negative convex functions belong to this class.

**Definition 3[5]:** A non $-$ negative function $f: [0, \infty) \to [0, \infty)$ is called $P-$convex (or that $f \in P(I)$) if for all $x, y \in [0, \infty)$ and $t \in [0, 1]$.
$$f(tx + (1-t)y) \leq f(x) + f(y)$$

**Definition 4[2, p. 288]:** A function $f: [0, \infty) \to \mathbb{R}$ is called an $s-$convex in second sense (or that $f \in K_s^2$), if for all $x, y \in I$ and $t \in (0, 1]$.
$$f(tx + (1-t)y) \leq t^s f(x) + (1-t)^s f(y)$$

Obviously, 1-convex function is convex.

**Definition 5[11]:** Let $I, J$ be intervals in $\mathbb{R}$, $(0,1) \subseteq J$ and let $h: J \subseteq \mathbb{R} \to \mathbb{R}$ be a non $-$ negative function, $h \not\equiv 0$. A non $-$ negative function $f: I \to \mathbb{R}$ is called $h-$convex $(or\ that\ f \in SX(h, I))$, if for all $x, y \in I$ and $t \in (0, 1)$
$$f(tx + (1-t)y) \leq h(t)f(x) + h(1-t)f(y) \tag{6}$$

If the inequality is reversed then $f$ is said to be $h$-concave and in this case $f$ belongs to the class $SV(h, I)$.

**Remark 1.**

i. If $h(t) = t$, then all the non – negative convex functions belong to the class $SX(h,I)$ and all non – negative concave functions belong to the class $SV(h,I)$.

ii. If $h(t) = 1$, then $SX(h,I) = Q(I)$

iii. If $h(t) = \frac{1}{t}$, then $SX(h,I) \supseteq P(I)$

iv. If $h(t) = t^s$, where $s \in (0,1)$, then $SX(h,I) \supseteq K_s^2$, respectively.

M. Z. Sarikaya et al. discussed in [8] some new inequalities of Hadamard type for $h$-convex functions. In [9], Sarikaya et. al. established a new Hadamard-type inequality for $h$-convex functions as:

**Theorem 5.** *Let $a, b \in I$ with $a < b$ and $f \in L[a,b]$. If $f \in SX(h,I)$ for $0 < t < 1$, then*

$$\frac{1}{2h\left(\frac{1}{2}\right)} f\left(\frac{a+b}{2}\right) \leq \frac{1}{b-a} \int_a^b f(x)\,dx \leq [f(a) + f(b)] \int_0^1 h(t)dt.$$

(7)

**Remark 2.** For $h(t) = t$ inequality (7) reduces to inequality (1) i:e.; Hermite-Hadamard inequality for convex function $f : [a,b] \to \mathbb{R}$.

In this article, we obtain new inequalities of Hermite-Hadamard type for functions whose absolute values of second derivatives are $h$-convex. We discuss refinements of $h$-Hermite-Hadamard inequality. Eventually, we have given some applications for special Means of real numbers.

**2. Main Results**

In [1], Barani et al. established inequalities for twice differentiable $P$-convex functions which are connected with Hadamards inequality, and they used the following lemma to prove their results:

**Lemma 1.** Let $f : I \subset \mathbb{R} \to \mathbb{R}$ be a twice differentiable function on $I^o$, $a, b \in I$ ($I^o$ is interior of $I$) with $a < b$: If $f'' \in L^1([a; b])$, then the following identity holds:

$$\frac{f(a) + f(b)}{2} - \frac{1}{b-a}\int_a^b f(x)dx$$

$$= \frac{(b-a)^2}{16}\int_0^1 (1-t^2)\left[f''\left(\frac{1+t}{2}a + \frac{1-t}{2}b\right) + f''\left(\frac{1-t}{2}a + \frac{1+t}{2}b\right)\right]dt.$$

In the following theorem, we will suggest some new upper bound for the right hand side of (1) for *h*-convex functions.

**Theorem 6.** Let $f : I \subset [0, \infty) \to \mathbb{R}$ be a twice differentiable function on $I^o$, such that $f'' \in L^1([a; b])$, where $a, b \in I^o$ with $a < b$: If $|f''|$ is *h*-convex on $[a, b]$ then for some fixed $t \in (0, 1)$ the following inequality holds:

$$\left|\frac{f(a) + f(b)}{2} - \frac{1}{b-a}\int_a^b f(x)dx\right|$$

$$\leq \frac{(b-a)^2}{16}\int_0^1 \left\{|1-t^2|\left(|f''(a)| + |f''(b)|\right) + |1 - (1-t^2)|\left|f''\left(\frac{a+b}{2}\right)\right|\right\}dt.$$

$$\leq \frac{(b-a)^2}{16}[|f''(a)| + |f''(b)|]\int_0^1 (1-t^2)\left\{h(t) + 2h\left(\frac{1}{2}\right)h(1-t)\right\}dt. \quad (8)$$

*Proof.* Since $|f''|$ is a *h*-convex function, by using Lemma 1 we get

$$\left|\frac{f(a) + f(b)}{2} - \frac{1}{b-a}\int_a^b f(x)dx\right|$$

$$\leq \frac{(b-a)^2}{16}\int_0^1 (1-t^2)\left\{\left|f''\left(\frac{1+t}{2}a + \frac{1-t}{2}b\right)\right| + \left|f''\left(\frac{1-t}{2}a + \frac{1+t}{2}b\right)\right|\right\}dt$$

$$\leq \frac{(b-a)^2}{16}\int_0^1 (1-t^2)\left\{h(t)(|f''(a)| + |f''(b)|) + 2h(1-t)\left|f''\left(\frac{a+b}{2}\right)\right|\right\}dt$$

$$\leq \frac{(b-a)^2}{16}[|f''(a)| + |f''(b)|]\int_0^1 (1-t^2)\left\{h(t) + 2h\left(\frac{1}{2}\right)h(1-t)\right\}dt.$$

A similar result is embodied in the following theorem:

**Theorem 7.** Let $f : I \subset [0, \infty) \to \mathbb{R}$ be a twice differentiable function on $I^o$, such that $f'' \in L^1([a; b])$, where $a, b \in I^o$ with $a < b$: If $|f''|^q$ for $q > 1$ with

$p = \frac{q}{q-1}$ is h-convex on $[a,b]$ then for some fixed $t \in (0,1)$ the following inequality holds:

$$\left|\frac{f(a)+f(b)}{2} - \frac{1}{b-a}\int_a^b f(x)dx\right| \leq \frac{(b-a)^2}{16.2^{\frac{1}{p}}}\beta^{\frac{1}{p}}\left(\frac{1}{2},p+1\right)\left(\int_0^1 h(t)dt\right)^{\frac{1}{q}} \times$$

$$\left[\left(|f''(a)|^q + \left|f''\left(\frac{a+b}{2}\right)\right|^q\right)^{\frac{1}{q}} + \left(|f''(b)|^q + \left|f''\left(\frac{a+b}{2}\right)\right|^q\right)^{\frac{1}{q}}\right]. \quad (9)$$

*Proof.* By using Lemma 1 and Hölder's inequality,

$$\left|\frac{f(a)+f(b)}{2} - \frac{1}{b-a}\int_a^b f(x)dx\right|$$

$$\leq \frac{(b-a)^2}{16}\int_0^1 (1-t^2)\left\{\left|f''\left(\frac{1+t}{2}a + \frac{1-t}{2}b\right)\right| + \left|f''\left(\frac{1-t}{2}a + \frac{1+t}{2}b\right)\right|\right\} dt$$

$$\leq \frac{(b-a)^2}{16}\left(\int_0^1 (1-t^2)^p dt\right)^{\frac{1}{p}} \times \left[\left(\int_0^1 \left|f''\left(\frac{1+t}{2}a + \frac{1-t}{2}b\right)\right|^q dt\right)^{\frac{1}{q}} + \right.$$

$$\left. \left(\int_0^1 \left|f''\left(\frac{1-t}{2}a + \frac{1+t}{2}b\right)\right|^q dt\right)^{\frac{1}{q}}\right].$$

Since $|f''|^q \in SX(h, I)$

$$\left|\frac{f(a)+f(b)}{2} - \frac{1}{b-a}\int_a^b f(x)dx\right| \leq \frac{(b-a)^2}{16}\left(\int_0^1 (1-t^2)^p dt\right)^{\frac{1}{p}} \times$$

$$\left[\left(\int_0^1 \left\{h(t)|f''(a)|^q + h(1-t)\left|f''\left(\frac{a+b}{2}\right)\right|^q\right\} dt\right)^{\frac{1}{q}} + \right.$$

$$\left. \left(\int_0^1 \left\{h(t)|f''(b)|^q + h(1-t)\left|f''\left(\frac{a+b}{2}\right)\right|^q\right\} dt\right)^{\frac{1}{q}}\right]$$

$$= \frac{(b-a)^2}{16}\left(\frac{1}{2}\beta\left(\frac{1}{2},p+1\right)\right)^{\frac{1}{p}}\left(\int_0^1 h(t)dt\right)^{\frac{1}{q}} \times$$

$$\left[\left(|f''(a)|^q + \left|f''\left(\frac{a+b}{2}\right)\right|^q\right)^{\frac{1}{q}} + \left(|f''(b)|^q + \left|f''\left(\frac{a+b}{2}\right)\right|^q\right)^{\frac{1}{q}}\right].$$

Here, Beta function is defined as:

$$\beta(u,v) = \int_0^1 t^{u-1}(1-t)^{v-1}dt, \quad u,v > 0.$$

And 
$$\int_0^1 h(t)dt = \int_0^1 h(1-t)dt \quad \blacksquare$$

Another similar result may be extended in the following theorem:

**Theorem 8.** Let $f : I \subset [0, \infty) \to \mathbb{R}$ be a twice differentiable function on $I^o$, such that $f'' \in L^1([a; b])$, where $a, b \in I^o$ with $a < b$: If $|f''|^q$ for $q \geq 1$ be an $h$-convex on $[a, b]$ then for some fixed $t \in (0, 1)$ the following inequality holds:

$$\left| \frac{f(a) + f(b)}{2} - \frac{1}{b-a} \int_a^b f(x)dx \right|$$

$$\leq \frac{(b-a)^2}{16} \left(\frac{2}{3}\right)^{\frac{1}{p}} \left[ \left(\int_0^1 \left\{ (1-t^2)|f''(a)|^q + t(2-t)\left|f''\left(\frac{a+b}{2}\right)\right|^q \right\} h(t)dt \right)^{\frac{1}{q}} \right.$$

$$\left. + \left(\int_0^1 \left\{ (1-t^2)|f''(b)|^q + t(2-t)\left|f''\left(\frac{a+b}{2}\right)\right|^q \right\} h(t)dt \right)^{\frac{1}{q}} \right] \quad (10)$$

*Proof:* By using Lemma 1 and well – known power – mean inequality,

$$\left| \frac{f(a)+f(b)}{2} - \frac{1}{b-a} \int_a^b f(x)dx \right|$$

$$\leq \frac{(b-a)^2}{16} \int_0^1 (1-t^2) \left\{ \left| f''\left(\frac{1+t}{2}a + \frac{1-t}{2}b\right) \right| + \left| f''\left(\frac{1-t}{2}a + \frac{1+t}{2}b\right) \right| \right\} dt$$

$$\leq \frac{(b-a)^2}{16} \left(\int_0^1 (1-t^2)dt\right)^{\frac{1}{p}} \left[ \left(\int_0^1 (1-t^2) \left| f''\left(\frac{1+t}{2}a + \frac{1-t}{2}b\right) \right|^q dt \right)^{\frac{1}{q}} + \right.$$

$$\left. \left(\int_0^1 (1-t^2) \left| f''\left(\frac{1-t}{2}a + \frac{1+t}{2}b\right) \right|^q dt \right)^{\frac{1}{q}} \right].$$

Since $|f''|^q \in SX(h, I)$

$$\left| \frac{f(a)+f(b)}{2} - \frac{1}{b-a} \int_a^b f(x)dx \right|$$

$$\leq \frac{(b-a)^2}{16} \left(\frac{2}{3}\right)^{\frac{1}{p}} \left[ \left(\int_0^1 (1-t^2) \left\{ h(t)|f''(a)|^q + h(1-t)\left|f''\left(\frac{a+b}{2}\right)\right|^q \right\} dt \right)^{\frac{1}{q}} + \right.$$

$$\left. \left(\int_0^1 (1-t^2) \left\{ h(t)|f''(b)|^q + h(1-t)\left|f''\left(\frac{a+b}{2}\right)\right|^q \right\} dt \right)^{\frac{1}{q}} \right]$$

$$= \frac{(b-a)^2}{16} \left(\frac{2}{3}\right)^{\frac{1}{p}} \left[ \left(\int_0^1 \left\{ (1-t^2)|f''(a)|^q + t(2-t)\left|f''\left(\frac{a+b}{2}\right)\right|^q \right\} h(t)dt \right)^{\frac{1}{q}} \right.$$

$$\left. \left(\int_0^1 \left\{ (1-t^2)|f''(b)|^q + t(2-t)\left|f''\left(\frac{a+b}{2}\right)\right|^q \right\} h(t)dt \right)^{\frac{1}{q}} \right].$$

which completes the proof. ∎

A parallel result for functions belonging to class of $h$-concave is stated as:

**Theorem 9.** Let $f : I \subset [0, \infty) \to \mathbb{R}$ be a twice differentiable function on $I^o$, such that $f'' \in L^1([a; b])$, where $a, b \in I^o$ with $a < b$: If $|f''|^q$ for $q \geq 1$ be an h-concave on $[a, b]$ then for some fixed $t \in (0, 1)$ the following inequality holds:

$$\left|\frac{f(a)+f(b)}{2} - \frac{1}{b-a}\int_a^b f(x)dx\right|$$

$$\leq \frac{(b-a)^2}{32}\left(\beta\left(\frac{1}{2}, p+1\right)\right)^{\frac{1}{p}}\left(\frac{1}{h\left(\frac{1}{2}\right)}\right)^{\frac{1}{q}}\left[\left|f''\left(\frac{a+3b}{4}\right)\right| + \left|f''\left(\frac{3a+b}{4}\right)\right|\right] \quad (11)$$

*Proof*: By employing Lemma 1 and Hölder's inequality,

$$\left|\frac{f(a)+f(b)}{2} - \frac{1}{b-a}\int_a^b f(x)dx\right|$$

$$\leq \frac{(b-a)^2}{16}\int_0^1 (1-t^2)\left\{\left|f''\left(\frac{1+t}{2}a + \frac{1-t}{2}b\right)\right| + \left|f''\left(\frac{1-t}{2}a + \frac{1+t}{2}b\right)\right|\right\} dt$$

$$\leq \frac{(b-a)^2}{16}\left(\int_0^1 (1-t^2)^p dt\right)^{\frac{1}{p}} \times \left[\left(\int_0^1 \left|f''\left(\frac{1+t}{2}a + \frac{1-t}{2}b\right)\right|^q dt\right)^{\frac{1}{q}} + \left(\int_0^1 \left|f''\left(\frac{1-t}{2}a + \frac{1+t}{2}b\right)\right|^q dt\right)^{\frac{1}{q}}\right].$$

Since $|f''|^q \in SV(h, I)$, therefore by inequality (7),

$$\int_0^1 \left|f''\left(\frac{1+t}{2}a + \frac{1-t}{2}b\right)\right|^q dt = \frac{1}{2h\left(\frac{1}{2}\right)}\left|f''\left(\frac{a+3b}{4}\right)\right|^q$$

and

$$\int_0^1 \left|f''\left(\frac{1-t}{2}a + \frac{1+t}{2}b\right)\right|^q dt = \frac{1}{2h\left(\frac{1}{2}\right)}\left|f''\left(\frac{3a+b}{4}\right)\right|^q,$$

which completes the proof. ∎

In order to prove our midpoint type inequalities, we need the following Lemma:

**Lemma 2.** Let $f : I \subset \mathbb{R} \to \mathbb{R}$ be twice differentiable function on $I°$, $a, b \in I°$ with $a < b$ and $f'' \in L^1([a, b])$, then the following identity holds:

$$\frac{1}{b-a}\int_a^b f(x)\, dx - f\left(\frac{a+b}{2}\right)$$

$$= \frac{(b-a)^2}{4}\int_0^1 m(t)\, [f''(t\,a + (1-t)b)dt + f''(t\,b + (1-t)a)]\, dt,$$

where $$m(t) = \begin{cases} t^2, & t \in \left[0, \frac{1}{2}\right) \\ (1-t)^2, & t \in \left[\frac{1}{2}, 1\right]. \end{cases}$$

**Theorem 10:** Let $f: I \subset [0, \infty) \to \mathbb{R}$ be a twice differentiable function on $I^\circ$, such that $f'' \in L^1([a, b])$, where $a, b \in I^\circ$ with $a < b$. If $|f''|$ is *h*-convex on $[a, b]$ for some fixed $t \in (0, 1)$ then the following identity holds:

$$\left| \frac{1}{b-a} \int_a^b f(x)dx - f\left(\frac{a+b}{2}\right) \right|$$
$$\leq (b-a)^2 \left[ \frac{|f''(a)| + |f''(b)|}{2} \right] \int_0^1 |m(t)| \, h(t) \, dt, \tag{12}$$

where $$\boldsymbol{m(t)} = \begin{cases} \boldsymbol{t^2}, & \boldsymbol{t \in \left[0, \frac{1}{2}\right)}, \\ \boldsymbol{(1-t)^2}, & \boldsymbol{t \in \left[\frac{1}{2}, 1\right]}. \end{cases}$$

**Theorem 11:** Let $f: I \subset [0, \infty) \to \mathbb{R}$ be a twice differentiable function on $I^\circ$, such that $f'' \in L^1([a, b])$, where $a, b \in I^\circ$ with $a < b$. If $|f''|^q$ for $q > 1$ with $p = \frac{q}{q-1}$ is *h*-convex on $[a, b]$ for some fixed $t \in (0, 1)$ then the following identity holds:

$$\left| \frac{1}{b-a} \int_a^b f(x)dx - f\left(\frac{a+b}{2}\right) \right|$$
$$\leq \frac{(b-a)^2}{8} \left( \frac{1}{2p+1} \right)^{\frac{1}{p}} \left( [|f''(a)|^q + |f''(b)|^q] \int_0^1 h(t)dt \right)^{\frac{1}{q}} \tag{13}$$

**Theorem 12:** Let $f: I \subset [0, \infty) \to \mathbb{R}$ be a twice differentiable function on $I^\circ$, such that $f'' \in L^1([a, b])$, where $a, b \in I^\circ$ with $a < b$. If $|f''|^q$ for $q \geq 1$ is *h*-convex on $[a, b]$ for some fixed $t \in (0, 1)$ then the following identity holds:

$$\left| \frac{1}{b-a} \int_a^b f(x)dx - f\left(\frac{a+b}{2}\right) \right|$$
$$\leq \frac{(b-a)^2}{2} \left( \frac{1}{12} \right)^{\frac{1}{p}} \left( [|f''(a)| + |f''(b)|] \int_0^1 |m(t)| \, h(t) \, dt \right)^{\frac{1}{q}}, \tag{14}$$

where
$$m(t) = \begin{cases} t^2, & t \in \left[0, \frac{1}{2}\right), \\ (1-t)^2, & t \in \left[\frac{1}{2}, 1\right]. \end{cases}$$

**Theorem 13:** let $f: I \subset [0, \infty) \to \mathbb{R}$ be a twice differentiable function on $I°$, such that $f'' \in L^1([a,b])$, where $a, b \in I°$ with $a < b$. if $|f''|^q$ for $p \geq 1$ is $h$-concave on $[a, b]$ for some fixed $t \in (0, 1)$ then the following identity holds:

$$\left| \frac{1}{b-a} \int_a^b f(x)dx - f\left(\frac{a+b}{2}\right) \right|$$

$$\leq \frac{(b-a)^2}{4} \left(\frac{1}{2p+1}\right)^{\frac{1}{p}} \left(\frac{1}{2h\left(\frac{1}{2}\right)}\right)^{\frac{1}{q}} \left|f''\left(\frac{a+b}{2}\right)\right|. \quad (15)$$

Proofs of theorems $10 - 13$ are bi-passed because the same are replica of the aforementioned theorems.

**Remark 3.** All of the above inequalities obviously hold for non-negative convex functions. Simply choose $h(t) = t$; $h(t) = \frac{1}{t}$; $h(t) = 1$ and $h(t) = t^s$ in each of those results, to get desired results respectively for non-negative convex functions, $Q(I)$; $P(I)$ and $k_s^2(I)$ : parallel development for $h$-concave functions may also be made likewise.

**Example 1.** For $h(t) = 1$; inequality (10) reduces inequality (3) and for $h(t) = t$; inequality (12) reduces inequality (2):

## 3. Refinement Of $h$-Hermite-Hadamard Inequality

To find new refined upper bound, integrating inequality (6) with respect to $t$,

$$\frac{1}{x-y} \int_{xa+(1-a)y}^{bx+(1-b)y} f(u)du \leq f(x) \int_a^b h(t)dt + f(y) \int_a^b h(1-t)dt.$$

For better right bound of hadamard inequality compare the above bound with usual one *i.e.*, $[f(a) + f(b)] \int_0^1 h(t)dt$. Let the above is less than the usual bound, that is,

$$f(x) \int_a^b h(t)dt + f(y) \int_a^b h(1-t)dt \leq [f(a) + f(b)] \int_0^1 h(t)dt$$

So, for the above being true, that is:

$$\int_a^b h(t)dt \leq \int_0^1 h(t)dt; \qquad \int_a^b h(1-t)dt \leq \int_0^1 h(t)dt,$$

taking $b = a + \lambda$ and using fact that $\int_0^1 h(t)dt = \int_0^1 h(1-t)dt$,

$$\int_a^{a+\lambda} h(t)dt \leq \int_0^1 h(t)dt; \qquad \int_a^{a+\lambda} h(1-t)dt \leq \int_0^1 h(1-t)dt. \quad (16)$$

From (16), it is cleared that $\lambda \leq 1$.

This means we have revamped the upper bound of Hermite Hadamard inequality for $h$-convex function, when the distance between $a$ and $b$ is almost one. The most interesting thing is that all linking work is with the interval $[0, 1]$ better than other. This discussion gives the following result.

**Theorem 14:** Let $f: [a + \lambda] \to \mathbb{R}$ be belonged to class $SX(h, I)$ for some $t \in (0, 1)$, $0 < \lambda \leq 1$ and $0 \leq a < 1$, then

$$\frac{1}{2h\left(\frac{1}{2}\right)} f\left(\frac{2a + \lambda}{2}\right) \leq \frac{1}{\lambda} \int_a^{a+\lambda} f(t)dt \leq \int_a^{a+\lambda} [f(x)h(t) + f(y)h(1-t)]dt \int_0^1 h(t)dt$$

An improvement of left inequality of (7) is given as:

**Theorem 15:** Let $f: I \subseteq \mathbb{R} \to \mathbb{R}$ be belonged to class $SX(h, I)$ and $a, b \in I°$ with $a < b$, then

$$\frac{1}{b-a} \int_a^b f(x)dx - \frac{1}{2h\left(\frac{1}{2}\right)} f\left(\frac{a+b}{2}\right)$$

$$\geq \left| \frac{1}{b-a} \int_a^b \left| \frac{f(x) + f(a+b-x)}{2} \right| dx - \frac{1}{2h\left(\frac{1}{2}\right)} \left| f\left(\frac{a+b}{2}\right) \right| \right| \geq 0.$$

*Proof:* By the $h$-convexity on $f$,

$$f\left(\frac{x+y}{2}\right) \le h\left(\frac{1}{2}\right)[f(x)+f(y)]$$

$$h\left(\frac{1}{2}\right)[f(x)+f(y)] - f\left(\frac{x+y}{2}\right) \ge \left|h\left(\frac{1}{2}\right)|f(x)+f(y)| - \left|f\left(\frac{x+y}{2}\right)\right|\right|$$

Let $x = t\,a + (1-t)b$ and $y = (1-t)a + t\,b$ for $t \in [0,1]$,

$$h\left(\frac{1}{2}\right)[f(t\,a + (1-t)b) + f((1-t)a + t\,b)] - f\left(\frac{a+b}{2}\right)$$

$$\ge \left|h\left(\frac{1}{2}\right)|f(t\,a + (1-t)b) + f((1-t)a + t\,b)| - \left|f\left(\frac{a+b}{2}\right)\right|\right|,$$

integrating from 0 to 1 over $t$,

$$h\left(\frac{1}{2}\right)\left[\int_0^1 f(t\,a + (1-t)b)\,dt + \int_0^1 f((1-t)a + t\,b)\,dt\right] - f\left(\frac{a+b}{2}\right)$$

$$\ge \left|h\left(\frac{1}{2}\right)\int_0^1 |f(t\,a + (1-t)b) + f((1-t)a + t\,b)|\,dt - \left|f\left(\frac{a+b}{2}\right)\right|\right|$$

$$\frac{2h\left(\frac{1}{2}\right)}{b-a}\int_a^b f(x)\,dx - f\left(\frac{a+b}{2}\right) \ge \left|\frac{2h\left(\frac{1}{2}\right)}{b-a}\int_a^b \left|\frac{f(x)+f(a+b-x)}{2}\right|\,dx - \left|f\left(\frac{a+b}{2}\right)\right|\right|,$$

which completes the proof. ∎

**Remark 4:** For $h(t) = t$, Theorem 15 coincides with Theorem 3

A refinement of second part of inequality (7) is given as:

**Theorem 16:** Let $f : I \subset \mathbb{R} \to \mathbb{R}$ be an *h*-convex function on the interval of real numbers $I$. Suppose $a; b \in I$ with $a < b$, then we have:

$$[f(a)+f(b)]\int_0^1 h(t)\,dt - \frac{1}{b-a}\int_a^b f(x)\,dx$$

$$\ge \begin{cases} 2f(a)\int_0^1 h(t)\,dt - \frac{1}{b-a}\int_a^b f(x)\,dx, & f(a) = f(b) \\ [f(a)+f(b)]\int_0^1 h(t)\,dt - \frac{1}{b-a}\int_a^b f(x)\,dx, & f(a) \ne f(b) \end{cases}$$

.

Proof: By the $h-$convexity on $f$ and continuity property of modulus

$f(t\,a + (1-t)b) \leq h(t)f(a) + h(1-t)f(b)$

$h(t)f(a) + h(1-t)f(b) - f(t\,a + (1-t)b)$
$$= |h(t)f(a) + h(1-t)f(b) - f(t\,a + (1-t)b)|$$
$$\geq ||h(t)f(a) + h(1-t)f(b)| - |f(t\,a + (1-t)b)|| \geq 0,$$

for all $a, b \in I$ and $t \in [0, 1]$:

integrating from 0 to 1 over $t$,

$f(a)\int_0^1 h(t)dt + f(b)\int_0^1 h(1-t)\,dt - \int_0^1 f(t\,a + (1-t)b)dt$
$$\geq \left|\int_0^1 |h(t)f(a) + h(1-t)f(b)|dt - \int_0^1 |f(t\,a + (1-t)b)|dt\right| \geq 0. \quad \blacksquare$$

**Remark 5:** For $h(t) = t$, Theorem **16** coincides with Theorem **4**.

**4. Applications to Some Special Means**

Let consider the applications to the following special means.

a). The arithmetic mean:

$$A = A(a, b) = \frac{a+b}{2}, \qquad a, b \in \mathbb{R},\ a, b > 0$$

b). The geometric mean:

$$G = G(a, b) = \sqrt{a\,b}, \qquad a, b \in \mathbb{R},\ a, b > 0$$

c). The harmonic mean:

$$H = H(a, b) = \frac{2\,a\,b}{a+b}, \qquad a, b \in \mathbb{R},\ a, b > 0$$

d). The logarithmic mean:

$$L(a,\ b) = \begin{cases} a & if\ a = b \\ \frac{b-a}{\ln b - \ln a}, & if\ a \neq b, \end{cases} \qquad a, b \in \mathbb{R},\ a, b > 0$$

e). The identric mean:

$$I(a,\ b) = \begin{cases} a & if\ a = b \\ \frac{1}{e}\left(\frac{b^b}{a^a}\right)^{\frac{1}{b-a}}, & if\ a \neq b, \end{cases} \qquad a, b \in \mathbb{R},\ a, b > 0$$

f). The generalized logarithmic mean:

$$L_n(a,\ b) = \begin{cases} a & if\ a = b \quad a, b \in \mathbb{R},\ a, b > 0, \\ \left[\frac{b^{n+1}-a^{n+1}}{(n+1)(b-a)}\right]^{\frac{1}{n}}, & if\ a \neq b, \quad n \in \mathbb{R} \backslash \{-1, 0\} \end{cases}$$

It is well known that $L_n$ is monotonically increasing over $n \in \mathbb{R}$ with $L_0 = I$ and $L_{-1} = L$. In particular, the following holds:

$$H \leq G \leq L \leq I \leq A.$$

**Proposition 1:** Let $q > 1$ with $q(p-1) = p$, $0 < a < b$, and $h(t) = t$, then

$$|L(a,b) - A(a,b)|$$

$$\leq \frac{(b-a)^2}{16} \beta^{\frac{1}{p}}\left(\frac{1}{2}, p+1\right) A\left(\left(a^q + G^q(a,b)\right)^{\frac{1}{q}}, \left(b^q + G^q(a,b)\right)^{\frac{1}{q}}\right).$$

***Proof***: The proof follows by Theorem 7 setting for $f(x) = e^x$. ∎

**Proposition 2:** Let $n \in \mathbb{N}$, $n \geq 2$ and $0 < a < b$, then

$$|\ A(a^n, b^n) - L_n^n(a,b)|$$

$$\leq \frac{(b-a)^2}{12}\left(\frac{1}{n(n+1)}\right)^{\frac{1}{q}} A\left(\left(\frac{3a^{n-2}+5A^{n-2}(a,b)}{8}\right)^{\frac{1}{q}}, \left(\frac{3b^{n-2}+5A^{n-2}(a,b)}{8}\right)^{\frac{1}{q}}\right).$$

*Proof*: Follows by Theorem 8, setting for $f(x) = x^n$ and $h(t) = t$. ∎

**Proposition 3:** Let $p > 1$ with $q(p-1) = p$ and $0 < a < b$, then

$$|L^{-1}(a^n, b^n) - A^{-1}(a,b)| \leq \frac{(b-a)^2}{4}\left(\frac{1}{2p+1}\right)^{\frac{1}{p}} H^{-\frac{1}{q}}(a^{3q}, b^{3q}).$$

*Proof*: Follows by Theorem 11 applied for $f(x) = \frac{1}{x}$ and $h(t) = t$. ∎

**Proposition 4:** Let $p > 1$ with $q(p-1) = p$ and $0 < a < b$, then

$$\left|\frac{I(a,b)}{A(a,b)}\right| \leq \exp\left[\frac{(b-a)^2}{4}\left(\frac{1}{2p+1}\right)^{\frac{1}{p}}\left(\frac{1}{2}\right)^{\frac{1}{q}} A^{-2}(a,b)\right].$$

***Proof***: Follows by Theorem 13 applying for $f(x) = -\ln x$ and $h(t) = t$. ∎